\documentstyle[12pt]{article}

\catcode`\@=11
\begin{document}
\centerline{\large\bf Universal ${\cal R}$-matrix  Of The Super Yangian Double $DY(gl(1|1))$  }
\vspace{0.8cm}
\centerline{\sf  $^a$Jin-fang Cai, $^{bc}$Shi-kun Wang,   
$^a$Ke Wu and $^a$Chi Xiong } 
\baselineskip=13pt
\vspace{0.5cm}
\centerline{$^a$ Institute of Theoretical Physics, Academia Sinica, }
\baselineskip=12pt
\centerline{ Beijing, 100080, P. R. China }
\vspace{0.3cm}
\centerline{ $^b$ CCAST (World Laboratory), P.O. Box 3730, Beijing, 100080, 
  P. R. China  }
\baselineskip=12pt
\vspace{0.3cm}
\centerline{ $^c$ Institute of Applied Mathematics, Academia Sinica,}
\baselineskip=12pt
\centerline{Beijing, 100080, P. R. China }
\vspace{0.9cm}
\begin{abstract} {Based on Drinfel$^{\prime}$d realization of super Yangian Double
 $DY\left(gl(1|1)\right)$, its pairing relations and universal
 ${\cal R}$-matrix are given. By taking evaluation representation of 
universal ${\cal R}$-matrix,
 another realization $L^{\pm}(u)$
of $DY\left(gl(1|1)\right)$ is obtained. 
These two realizations of $DY\left(gl(1|1)\right)$
are related by the supersymmetric extension of Ding-Frenkel map.  }
\end{abstract}
\vspace{2cm}

Yangian algebra was introduced by Drinfel$^{\prime}$d\cite{DRI1,DRI2}.
The quantum 
double of Yangian consists of Yangian itself and its dual with opposite 
comultiplication. There are three methods to define  the Yangian and Yangian 
double: Drinfel$^{\prime}$d-Jmbo \cite{DRI1,JIM}, Drinfel$^{\prime}$d new realization \cite{DRI2} 
and RS 
 approach \cite{RS} (or FRT approach \cite{FRT}
in the case of without center extension). 
 The explicit isomorphism between  Drinfel$^{\prime}$d new realization 
 and RS realization
of Yangian double can be 
 established through Gauss decomposition, the similar method used by Ding 
 and Frenkel in the discussions of quantum Affine algebra. \cite{DF}. 
The property of Yangian double, such as quasi-triangular properties
 and equivalence of Drinfel$^{\prime}$d and RS realization  
 was studied well in some papers \cite{KT1,KT2, IOH}.  
Although the Drinfel$^{\prime}$d realization of super Yangian double \cite{CAI, YZZ1} 
was constructed 
by means of RS method and Gauss decomposition, 
the quasi-triangular property, such as  universal ${\cal R}$-matrix of super 
Yangian double  has not been studied yet. In this paper, we find the 
Hopf pairing relations between super Yangian $Y\left(gl(1|1)\right)$
and its dual, then construct
the universal ${\cal R}$-matrix of $DY\left(gl(1|1)\right)$. By taking evaluation
representation, we get the FRT realization of $DY\left(gl(1|1)\right)$. 

\vspace{1cm}

Super Yangian double $DY\left(gl(1|1)\right)$ is the Hopf algebra generated
by elements\\
 (Drinfel$^{\prime}$d generators)
 $e_n, f_n, h_n, k_n,  n\in {\bf Z} $ which satisfy the following
multiplication relations
\begin{eqnarray}
&&[h_m ~,~ h_n]=[h_m ~,~ k_n]=[k_m ~,~k_n]=0 \nonumber \\
&&[k_m ~,~ e_n]=[k_m ~,~ f_n] =0 \nonumber \\
&&[h_0~,~e_n]=-2e_n ~,~ [k_0~,~f_n]=2f_n \nonumber \\
&&[h_{m+1}~,~e_n]-[h_m~,~e_{n+1}]+\{h_m~,~e_n\}=0 \label{dg}\\
&&[h_{m+1}~,~f_n]-[h_m~,~f_{n+1}]-\{h_m~,~f_n\}=0 \nonumber \\
&&\{e_m~,~e_n\}=\{f_m~,~f_n\}=0 \nonumber \\
&&\{e_m~,~ f_n\}=-k_{m+n} \nonumber 
\end{eqnarray}
They could also be written as generating functions 
(or Drinfel$^{\prime}$d currents)
\begin{eqnarray}
&&E^{\pm}(u)=\pm \sum_{n \ge 0 \atop n<0} e_n u^{-n-1} ~,~~~~~~~~
F^{\pm}(u)=\pm \sum_{n \ge 0 \atop n<0} f_n u^{-n-1} \\
&&H^{\pm}(u)=1\pm \sum_{n \ge 0 \atop n<0} h_n u^{-n-1} ~,~~~~
K^{\pm}(u)=1\pm \sum_{n \ge 0 \atop n<0} k_n u^{-n-1}
\end{eqnarray}
and
\begin{eqnarray}
&&E(u)=E^+(u)-E^-(u) ~,~~~~~ F(u)=F^+(u)-F^-(u)
\end{eqnarray}
then the relations (1) look as follows
\begin{eqnarray}
&&[H^{\sigma}(u)~,~H^{\rho}(v)]=[H^{\sigma}(u)~,~K^{\rho}(v)]
=[K^{\sigma}(u)~,~K^{\rho}(v)]=0 ,~~~~\forall \sigma,\rho=+,- \nonumber \\
&&[K^{\pm}(u)~,~E(v)]=[K^{\pm}(u)~,~F(v)]=0\nonumber  \\
&&\{E(u)~,~E(v)\}=\{F(u)~,~F(v)\}=0 \nonumber \\
&&H^{\pm}(u)E(v)=\frac{u-v-1}{u-v+1}E(v)H^{\pm}(u) \label{dc}  \\
&&H^{\pm}(u)F(v)=\frac{u-v+1}{u-v-1}F(v)H^{\pm}(u) \nonumber  \\
&&\{E(u)~,~F(v)\}=\delta(u-v)[K^-(v)-K^+(u)]\nonumber 
\end{eqnarray}
in which $\delta(u-v)=\sum_{k\in Z}u^kv^{-k-1}$. The comultiplication structure
for $DY\left(gl(1|1)\right)$ is given by
\begin{eqnarray}
&&\triangle \left(E^{\pm}(u)\right)= E^{\pm}(u)\otimes 1+ 
   H^{\pm}(u) \otimes E^{\pm}(u) \nonumber \\
&&\triangle \left(F^{\pm}(u)\right)=1\otimes F^{\pm}(u)+
F^{\pm}(u)\otimes H^{\pm}(u) \nonumber \\
&&\triangle \left(K^{\pm}(u)\right)=K^{\pm}(u)\otimes K^{\pm}(u) \label{dco} \\
&&\triangle \left(H^{\pm}(u)\right)= H^{\pm}(u)\otimes H^{\pm}(u)-
 2F^{\pm}(u-1) H^{\pm}(u)\otimes H^{\pm}(u)E^{\pm}(u-1)) \nonumber 
\end{eqnarray}

As a quantum double, $DY\left(gl(1|1)\right)$ consists of the super 
Yangian $Y\left(gl(1|1)\right)$  and its dual $Y^*\left(gl(1|1)\right)$ 
with opposite comultiplication. The super Yangian $Y\left(gl(1|1)\right)$ 
is generated by $E^+(u), F^+(u)$,\\ 
$ H^+(u), K^+(u)$ and 
$Y^*\left(gl(1|1)\right)$ is generated by $E^-(u), F^-(u), H^-(u), K^-(u)$.
There exists a Hopf pairing relation between $Y\left(gl(1|1)\right)$ and 
$Y^*\left(gl(1|1)\right)$ : $<~~,~~>$ which satisfies the conditions
\begin{eqnarray}
&&<ab~,~ c^*d^*>=<\triangle(ab)~,~c^*\otimes d^*>=<b\otimes a ~,~\triangle
(c^*d^*)>
\end{eqnarray}
for any $a,b \in Y\left(gl(1|1)\right) $ and $c^*, d^* \in 
Y^*\left(gl(1|1)\right)$. We find that this pairing relation can be 
written as
\begin{eqnarray}
&&<E^+(u)~,~F^-(v)>=\frac{1}{u-v},~~~~~~~~<F^+(u)~,~E^-(v)>=\frac{1}{u-v} \\
&&<H^+(u)~,~K^-(v)>=\frac{u-v-1}{u-v+1},~~<K^+(u)~,~H^-(v)>=\frac{u-v-1}{u-v+1}
\end{eqnarray}

As the same discussion for $DY(sl_2)$ \cite{KT1}, the universal 
${\cal R}$-matrix
for $DY\left(gl(1|1)\right)$ has the following form
\begin{eqnarray}
&&{\cal R}={\cal R}_+{\cal R}_1{\cal R}_2{\cal R}_-  \label{rmatrix}
\end{eqnarray}
where
\begin{eqnarray}
&&{\cal R}_+=\prod_{n\ge 0}^{\rightarrow}\exp(- e_n\otimes f_{-n-1})  \\
&&{\cal R}_-=\prod_{n\ge 0}^{\leftarrow}\exp(- f_n\otimes e_{-n-1})  \\
&&{\cal R}_1=\prod_{n\ge 0} \exp \left\{ {\rm Res}_{u=v}\left[(-1)
\frac{\rm d}{{\rm d}u}({\rm ln}H^+(u))\otimes 
{\rm ln}K^-(v+2n+1)\right]\right\} \\
&&{\cal R}_{2}=\prod_{n\ge 0} \exp \left\{ {\rm Res}_{u=v}\left[(-1)
\frac{\rm d}{{\rm d}u}({\rm ln}K^+(u))\otimes 
{\rm ln}H^-(v+2n+1)\right]\right\} 
\end{eqnarray}
here we have used the notations
\begin{eqnarray}
&&{\rm Res}_{u=v}\left(A(u)\otimes B(v)\right)=\sum_k a_k\otimes b_{-k-1}
\end{eqnarray}
for $A(u)=\sum_k a_k u^{-k-1}$ and $B(u)=\sum_k b_k u^{-k-1}$. 

From the quasi-triangular property of the double, the universal ${\cal R}$-matrix
satisfies
\begin{eqnarray}
&&{\cal R}_{12}\cdot {\cal R}_{13}\cdot {\cal R}_{23}={\cal R}_{23}\cdot {\cal R}_{13}\cdot {\cal R}_{12} \label{rrr} \\
\vspace{3mm}
&&(\triangle \otimes id){\cal R}={\cal R}_{13} \cdot {\cal R}_{23},~~~~~
(id \otimes \triangle ){\cal R}={\cal R}_{13} \cdot {\cal R}_{12} \label{rco}
\end{eqnarray}
 In dealing with the tensor product in the graded case, we must use the form
$(A\otimes B)\cdot (C\otimes D)=(-1)^{P(B)P(c)} AC\otimes BD$, $P(B)=0, 1 $ for 
$B$ is bosonic and fermionic respectively. 

Let $\rho _x$ be taking  two-dimensional  evaluation representation for   
$DY\left(gl(1|1)\right)$:
\begin{eqnarray}
&&\rho _x(e_n) =\left(\begin{array}{ll}0& 0 \\x^n&0 \end{array}\right),~~~~~~~~~~~~~~
\rho _x(f_n) =\left(\begin{array}{ll}0& x^n \\0&0 \end{array}\right)  \\
&&\rho _x(h_n) =\left(\begin{array}{ll}x^n&0 \\0&-x^n \end{array}\right),~~~~~~~~~
\rho _x(k_n) =\left(\begin{array}{ll}-x^n& 0 \\0&-x^n \end{array}\right) 
\end{eqnarray}
and let
\begin{eqnarray}
&&L^+(x)=(\rho_x\otimes id)({\cal R}^{21})^{-1}, ~~~~~~ L^-(x)=(\rho_x\otimes id){\cal R} \\
&&R^+(x-y)=(\rho_x\otimes \rho_y)({\cal R}^{21})^{-1}, ~~~~
R^-(x-y)=(\rho_x\otimes \rho_y){\cal R}
\end{eqnarray}
then from (\ref{rmatrix}), we have 
\begin{eqnarray}
&&L^+(x)=\left(\begin{array}{lr}1&0 \\F^+(x)&1\end{array}\right)
    \left(\begin{array}{lr}k_1^+(x)&0\\0&k_2^+(x)\end{array}\right)
    \left(\begin{array}{lr}1&E^+(x)\\0&1 \end{array}\right) \label{lp} \\
&&L^-(x)=\left(\begin{array}{lr}1&0 \\F^-(x)&1\end{array}\right)
    \left(\begin{array}{lr}k_1^-(x)&0\\0&k_2^-(x)\end{array}\right)
    \left(\begin{array}{lr}1&E^-(x)\\0&1 \end{array}\right) \label{lm}
\end{eqnarray}
here
\begin{eqnarray}
&&k_1^+(x)=\prod_{n\ge 0}\frac{K^+(x-2n-2)}{K^+(x-2n-1)}\frac{H^+(x-2n)}
{H^+(x-2n-1)}  \\ 
&&k_2^+(x)=\prod_{n\ge 0}\frac{K^+(x-2n)}{K^+(x-2n-1)}\frac{H^+(x-2n)}
{H^+(x-2n-1)}  \\
&&k_1^-(x)=\prod_{n\ge 0}\frac{K^+(x+2n+1)}{K^+(x+2n)}\frac{H^+(x+2n+1)}
{H^+(x+2n+2)}  \\
&&k_2^-(x)=\prod_{n\ge 0}\frac{K^+(x+2n+1)}{K^+(x+2n+2)}\frac{H^+(x+2n+1)}
{H^+(x+2n+2)}
\end{eqnarray}
and
\begin{eqnarray}
&&R^{\pm}(x-y)=\rho^{\pm}(x-y)\left( \begin{array}{lccr}1&0&0&0\\
0&\frac{x-y}{x-y+1}&\frac{1}{x-y+1}&0\\0&\frac{1}{x-y+1}&\frac{x-y}{x-y+1}&0\\
0&0&0&\frac{x-y-1}{x-y+1} \end{array}\right)
\end{eqnarray}
here
\begin{eqnarray}
&&\rho^+(x)=\prod_{n\ge 0}\frac{(x-2n-3)(x-2n-1)^2(x-2n+1)}
{(x-2n-2)^2(x-2n)^2}  \\
&&\rho^-(x)=\prod_{n\ge 0}\frac{(x+2n)^2(x+2n+2)^2}{(x+2n-1)(x+2n+1)^2
(x+2n+3)}
\end{eqnarray}

From (\ref{rrr}), the relations among $R^{\pm}(x-y)$ and $L^{\pm}(x)$ can be obtained 
\begin{eqnarray}
&&R_{ij,ab}(x-y)R_{ak,pc}(x-z)R_{bc,qr}(y-z)
(-1)^{(P(a)-P(p))P(b)}   \nonumber    \\
\vspace{3mm}
&&{\hspace{2cm}}=(-1)^{P(e)(P(f)-P(r))}
R_{jk,ef}(y-z)R_{if,dr}(x-z)R_{de,pq}(x-y) \\
\vspace{3mm}
&&R^{\pm}_{ij,mn}(u-v)L^{\pm}_{mk}(u)L^{\pm}_{nl}(v)
(-1)^{P(k)(P(n)+P(l))}\nonumber \\
\vspace{3mm}
&& \hspace{2cm}
=
(-1)^{P(i)(P(j)+P(q))} L^{\pm}_{jq}(v) L^{\pm}_{ip}(u)R^{\pm}_{pq,kl}(u-v)
  \label{rll1}\\
  \vspace{3mm}
&&R^{-}_{ij,mn}(u-v)L^-_{mk}(u) 
L^+_{nl}(v)(-1)^{P(k)(P(n)+P(l))}\nonumber \\
\vspace{3mm}
&&\hspace{2cm} =
(-1)^{P(i)(P(j)+P(q))} L^+_{jq}(v) L^-_{ip}(u)R^-_{pq,kl}(u-v)
\label{rll2}
\end{eqnarray}
The comultiplication structure for current $L_{ij}^{\pm}(u)$ are got from
(\ref{rco})
\begin{eqnarray}
&&\triangle\left(L_{ij}^{\pm}(u)\right)=\sum_{k=1,2}(-1)^{(i+k)(k+j)}
L_{kj}^{\pm}(u)\otimes L_{ik}^{\pm}(u) \label{lco}
\end{eqnarray}

The relations (\ref{rll1}, \ref{rll2}) and (\ref{lco}) is another defining
of super Yangian double, which is usually referred to super version of 
FRT \cite{FRT} construction method. If we start from (\ref{rll1}, \ref{rll2})
 and (\ref{lco}) to define super Yangian double, using the decomposition
(\ref{lp}), (\ref{lm}) and setting $K^{\pm}(u)=k_1^{\pm}(u)^{-1}k_2^{\pm}(u),
H^{\pm}(u)=k_1^{\pm}(u)k_2^{\pm}(u-1)$, we can also rediscover the 
Drinfel$^{\prime}$d's currents or generators realization of the super 
Yangian double (\ref{dg}),(\ref{dc}) and (\ref{dco}) \cite{CAI} .

\end{document}